%% file: mt-paper.tex
\documentclass[10pt]{article}

\usepackage[a4paper]{geometry}
\usepackage{indentfirst}
\usepackage{hyperref}
\usepackage{amssymb, amsmath, amsthm}
\usepackage{graphicx}
\usepackage[labelsep=none]{caption}
\usepackage[labelformat=simple]{subcaption}

\hypersetup{bookmarksnumbered, pdfstartview=FitH, pdfpagemode=UseNone}

\newcommand{\llbr}{[\![}
\newcommand{\rrbr}{]\!]}
\newcommand{\within}{\mathbin{\upharpoonright}}

\newcommand{\mynewtheorem}[2]{\newtheorem{#1}{\indent #2}}
\newcommand{\myalttheorem}[2]{\newtheorem*{#1}{\indent #2}}
\mynewtheorem{conjecture}{Conjecture}
\myalttheorem{conjecture*}{Conjecture}
\mynewtheorem{lemma}{Lemma}
\mynewtheorem{proposition}{Proposition}
\myalttheorem{proposition*}{Proposition}
\mynewtheorem{theorem}{Theorem}
\myalttheorem{theorem*}{Theorem}
\newenvironment{myproof}[1][Proof]{\begin{proof}[\indent #1]}{\end{proof}}

\begin{document}

\title{\textbf{Minimum-Turn Tours of Even Polyominoes}}
\author{Nikolai Beluhov}
\date{}

\maketitle

\begin{abstract} Let $P$ be a connected bounded region in the plane formed out of $2 \times 2$ blocks joined by their sides. Peng and Rascoussier conjectured that all minimum-turn Hamiltonian cycles of $P$ exhibit a certain regular structure. We prove this conjecture in the special case when $P$ is a topological disk. The proof proceeds in two phases -- a ``downward'' phase where we break apart an irregular Hamiltonian cycle into a collection of shorter cycles; and an ``upward'' phase where we put it back together in a different way so that, overall, the number of turns in it decreases. \end{abstract}

\input{mt-paper-01-intro}

\input{mt-paper-02-prelim}

\input{mt-paper-03-tg}

\input{mt-paper-04-down}

\input{mt-paper-05-up}

\section*{Acknowledgements}

The author is thankful to Xiao Peng and Florian Rascoussier for telling him about their conjecture.

The present paper was written in the course of the author's PhD studies under the supervision of Professor Imre Leader. The author is thankful also to Prof.\ Leader for his unwavering support.

\end{document}

%% file: mt-paper-01-intro.tex
\section{Introduction} \label{intro}

Consider an infinite grid of unit square cells. Let $S$ be a connected bounded region in this grid, formed out of finitely many cells joined together by their sides. The \emph{adjacency graph} of $S$ is the graph whose vertices are the centres of the cells of $S$, and where two centres are joined by an edge if the corresponding cells are neighbours by side. A \emph{tour} of $S$ is a Hamiltonian cycle in the adjacency graph of $S$.

We write $2S$ for the region obtained from $S$ when it is scaled up by a factor of $2$; or, informally, when each cell of $S$ is replaced with a $2 \times 2$ block of four cells. Let $P = 2S$. We call regions $P$ of this form \emph{even}.

\begin{figure}[ht] \null \hfill \begin{subfigure}[c]{100pt} \centering \includegraphics{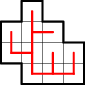} \caption{} \end{subfigure} \hfill \begin{subfigure}[c]{100pt} \centering \includegraphics{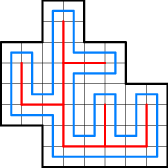} \caption{} \end{subfigure} \hfill \null \caption{} \label{tt} \end{figure}

It is well-known that an even region always admits a tour. \cite{GR} In fact, each spanning tree in the adjacency graph of $S$ maps onto a tour of $P$. Figure \ref{tt} shows one example, with a spanning tree of $S$ in red on the left and the associated tour of $P$ in blue on the right. We spell out the full description of the mapping in Section \ref{prelim}. We call the tours of $P$ obtained in this way \emph{regular}.

Suppose now that we wish to find a tour of $P$ which makes as few turns as possible. This is a natural quantity to optimise for in the original setting of \cite{GR} -- i.e., the setting of a robot surveying a connected bounded region in the plane. It might be that the robot has to slow down in order to make a turn, and so a tour with fewer turns would ensure that the survey is completed more promptly. Or perhaps the robot's memory is limited, and so we prefer tours which require fewer instructions of the form ``go this many steps in that direction'' to specify.

Even outside of this context, in pure mathematics, minimum-turn tours have attracted some interest because they arise naturally as auxiliary structures in other combinatorial problems. For example, both of \cite{TZ} and \cite{B} make use of the minimum-turn tours of rectangular regions.

The regular tours of $P$ are readily available to us because we can manufacture them ``cheaply'' out of the spanning trees of $S$. So it might be tempting to limit our search to them. However, it is far from obvious whether the smallest possible number of turns is always attained by some regular tour. Conceivably, it might be that all optimal tours are irregular. Indeed, it seems intuitively plausible that for many regions $S$ the regular tours of $P$ are a tiny subset of the set of all tours~of~$P$.

We arrive at the following question: Does $P$ always admit a regular tour which attains the smallest possible number of turns? Peng and Rascoussier \cite{PR} conjectured recently that, in fact, a much stronger claim is true:

\begin{conjecture*} \textnormal{(Peng and Rascoussier)} Let $P$ be an even region. Then all minimum-turn tours of $P$ are regular. \end{conjecture*}

Our goal will be to prove this conjecture in the special case when $P$ is a topological disk; or, informally, when $P$ is free of holes. We call such regions \emph{polyominoes}. We note that, even in this tamer setting, one already runs into multiple difficulties because the absence of holes does not prevent the boundary of $P$ from being ``arbitrarily complicated''.

One generalisation of the notion of a tour is going to be of particular importance to us. We define a \emph{pseudotour} of $P$ to be a cycle decomposition of the adjacency graph of $P$; i.e., a collection of pairwise disjoint cycles in this graph which together form a partitioning of the cells of $P$. Of course, a tour is just a pseudotour which consists of a single cycle.

We can now give an outline of our strategy. We begin with any irregular tour of $P$. By iterating a certain procedure, step by step we transform this tour into a pseudotour where the number of turns has been made smaller at the cost of the number of cycles increasing. Then, by iterating a different procedure, step by step we transform this pseudotour back into a tour -- but this time around we make the number of cycles smaller at the cost of allowing the number of turns to increase once again. We tally everything up, and we find that by going first ``down'' and then ``up'' in this manner, we have managed to save at least a couple of turns overall. Hence, an irregular tour of $P$ can never be optimal.

The rest of the paper is structured as follows: Section \ref{prelim} covers some preliminaries. Section~\ref{tg} defines one natural auxiliary structure which will be crucial to our argument. Then, in Section~\ref{down}, we describe the ``downward'' procedure and we study its behaviour. Finally, in Section \ref{up}, we describe the ``upward'' procedure and we complete the proof.

%% file: mt-paper-02-prelim.tex
\section{Preliminaries} \label{prelim}

Let $\mathcal{G}$ be our infinite grid of unit square cells. We formalise $\mathcal{G}$ as the graph whose vertices are all integer points in the plane and whose edges are all unit segments joining these points. We will also be working with the dual grid $\mathcal{H}$, formalised as the graph whose vertices are the centres of the cells of $\mathcal{G}$ and whose edges are all unit segments joining these centres. Of course, dually, the vertices of $\mathcal{G}$ coincide with the centres of the cells of $\mathcal{H}$.

Our reasoning will involve both graphs and polygons. The term ``vertex'' exists in both of these contexts, but with different meanings. To avoid confusion, from now on we will use ``vertex'' solely in its combinatorial sense, as ``vertex of a graph''; and we will use ``corner'' to mean ``vertex of a polygon''.

Let $S$ be a polyomino which lives in $\mathcal{G}$ -- i.e., which is the union of several cells of $\mathcal{G}$. We write $\partial S$ for the boundary of $S$; this will always be a cycle in $\mathcal{G}$. Conversely, if $C$ is a cycle in $\mathcal{G}$, we write $[C]$ for the polyomino enclosed by $C$. The same notations apply to $\mathcal{H}$ as well.

Let $a_1$, $a_2$, $\ldots$, $a_k$ be the corners of $S$, in this order on its boundary. Then we let $P = 2S$ be the polyomino with corners, in this order on its boundary, $2a_1$, $2a_2$, $\ldots$, $2a_k$. Of course, as defined, $P$ will also live in $\mathcal{G}$.

We introduce doubling for graphs, too. Let $H$ be a subgraph of $\mathcal{H}$. We write $2H$ for the subgraph of $\mathcal{G}$ obtained when we scale up $H$ by a factor of $2$. Or, more formally: For each edge $u$---$v$ of $H$, we take the path $2u$---$(u + v)$---$2v$ in $\mathcal{G}$; and then we form $2H$ out of all such paths. Given a subgraph $G$ of $\mathcal{G}$, we write $\llbr G \rrbr$ for the union of the cells of $\mathcal{H}$ centred at the vertices of~$G$.

We can now describe the mapping which transforms the spanning trees of $S$ into tours of $P$. Let $T$ be a spanning tree of $S$; so, in particular, $T$ is a subgraph of $\mathcal{H}$. It is intuitively obvious from Figure \ref{tt} what is going on -- the tour ``goes around'' $2T$. However, for our purposes it will be useful to spell out the details of the construction a bit more rigorously. Since $T$ is a tree, the region $\llbr 2T \rrbr$ is a polyomino. We claim that $\partial \llbr 2T \rrbr$ is a tour of $P$.

Consider the centre $p$ of any cell of $P$. It suffices to check that $p$ lies on $\partial \llbr 2T \rrbr$. Call a cell of $\mathcal{H}$ \emph{even} when both coordinates of its centre are even; and \emph{odd} when both coordinates of its centre are odd. Since $T$ spans $S$, we get that all odd cells of $\mathcal{H}$ contained within $P$ are also contained within $\llbr 2T \rrbr$. On the other hand, by the definition of $\llbr 2T \rrbr$ we get that it does not contain any even cells of $\mathcal{H}$.

Let $c'$ be the unique odd cell of $\mathcal{H}$ such that $p$ is a corner of it, and let $c''$ be the unique even one. Since $P$ is even, $c'$ is contained within $P$, and hence also within $\llbr 2T \rrbr$; but $c''$ is outside of~$\llbr 2T \rrbr$. We conclude that $p$ does indeed lie on $\partial \llbr 2T \rrbr$ -- thereby confirming our claim.

%% file: mt-paper-03-tg.tex
\section{The Turn Graph} \label{tg}

Let $\mathcal{U}$ be the graph whose vertices are all integer points in the plane and whose edges are all segments of length $\sqrt{2}$ joining these points; i.e., the graph formed by the vertices of $\mathcal{G}$ and the diagonals of the cells of $\mathcal{G}$. We write $\mathcal{U} \within P$ for the subgraph of $\mathcal{U}$ formed by the corners and the diagonals of the cells of $P$; we think of this as the ``restriction'' of $\mathcal{U}$ to $P$.

Let $E$ be a pseudotour of $P$. We define the \emph{turn graph} $U$ of $E$ to be the subgraph of $\mathcal{U} \within P$ formed by those diagonals in the cells of $P$ which ``bisect'' the turns of $E$. Or, more formally: For each turn $u'$---$v$---$u''$ of $E$, we take the diagonal joining $w' = (u' + u'')/2$ and $w'' = 2v - w'$; and then we form $U$ out of all such diagonals. Observe that, crucially, the number of turns in $E$ equals the number of edges in $U$.

\begin{figure}[ht] \null \hfill \begin{subfigure}[c]{100pt} \centering \includegraphics{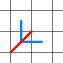} \caption{} \label{dia} \end{subfigure} \hfill \begin{subfigure}[c]{100pt} \centering \includegraphics{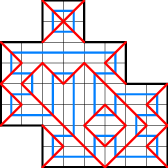} \caption{} \label{turn} \end{subfigure} \hfill \null \caption{} \label{dia-turn} \end{figure}

For example, Figure \ref{dia-turn} shows one turn in isolation on the left, alongside a pseudotour on the right together with its turn graph. In both instances, the edges of the pseudotour are coloured in blue and the edges of the associated turn graph are coloured in red.

Since each cell of $P$ contains at most one turn of $E$, we get that two diagonals in the same cell can never belong to $U$ simultaneously. Hence, $U$ is free of self-intersections. Our next observation regarding $U$ is only slightly more subtle:

\begin{lemma} \label{deg} Each corner of $P$ is an odd-degree vertex of $U$; and the other vertices of $U$ are all of even degree in $U$. \end{lemma} 

\begin{myproof} Consider any vertex $v$ of $\mathcal{G}$ in the interior of $P$. We must show that $v$ is incident with an even number of edges of $U$. (Possibly zero.) This is straightforward to do by going through all possible cases for the $12$ edges of $\mathcal{H}$ incident with $v + (\pm 1/2, \pm 1/2)$ being part of $E$ or not; but there is also a different argument which avoids case exhaustion altogether.

Let $c$ be the cell of $\mathcal{H}$ centred at $v$. Then $v$ is incident with an edge of $U$ which passes through the top right corner of $c$ if and only if the top and right sides of $c$ are in the same state relative to $E$; i.e., either both belonging to $E$ or both not belonging to $E$. Similar reasoning applies to the other three corners of $c$ as well; and, as we go around $c$ once, the states of its sides relative to $E$ must change an even number of times.

The cases when $v$ is either a corner of $P$ or otherwise on the boundary of $P$ are both handled analogously. \end{myproof}

This means that we can view $U$ as the union (edge-disjoint but not necessarily vertex-disjoint) of several paths and cycles, with the paths connecting the corners of $P$ by pairs.

We proceed to show that the properties of $U$ which we have worked out thus far in fact characterise the turn graphs completely.

\begin{lemma} \label{tgp} Consider any subgraph $U$ of $\mathcal{U} \within P$ such that: (i) $U$ is free of self-intersections; (ii) Each corner of $P$ is an odd-degree vertex of $U$; and (iii) The other vertices of $U$ are all of even degree in $U$. Then $U$ is the turn graph of some pseudotour $E$ of $P$. \end{lemma} 

\begin{myproof} By condition (i), the union of $\partial P$ with $U$ is a planar graph. (Or, strictly speaking, a planar graph embedding.) Its bounded faces form a partitioning of $P$ into sub-regions. By condition (iii), we can colour these sub-regions in black and white so that two sub-regions which neighbour one another across an edge of $U$ are always of opposite colours. By condition (ii), furthermore all horizontal sides of $P$ will belong to sub-regions of the same colour -- say, white -- and all vertical sides of $P$ will belong to sub-regions of the opposite colour -- i.e., black. We can now form $E$ out of all horizontal edges of $\mathcal{H}$ contained within the white sub-regions and all vertical edges of $\mathcal{H}$ contained within the black sub-regions. \end{myproof}

Or, in other words, there is a one-to-one correspondence between the pseudotours of $P$ and the subgraphs of $\mathcal{U} \within P$ which satisfy conditions (i)--(iii) of Lemma \ref{tgp}.

Clearly, $\mathcal{U}$ consists of two connected components. One of them spans all integer points $(x, y)$ with $x + y$ even, and we denote it by $\mathcal{U}_\text{Even}$; similarly, the other one spans all integer points $(x, y)$ with $x + y$ odd, and we denote it by $\mathcal{U}_\text{Odd}$. We define $\mathcal{U}_\text{Even} \within P$ and $\mathcal{U}_\text{Odd} \within P$ analogously to $\mathcal{U} \within P$.

We are ready now to formulate one handy criterion for regularity in terms of the turn graph.

\begin{lemma} \label{even} A tour $E$ of $P$ is regular if and only if its turn graph $U$ is a subgraph of $\mathcal{U}_\textnormal{Even}$. \end{lemma}

\begin{myproof} Consider the centre $p$ of any cell of $P$. Then, as in Section \ref{prelim}, we get that $p$ is a corner of exactly one odd cell $c'$ of $\mathcal{H}$ and exactly one even cell $c''$ of $\mathcal{H}$. Define a cell of $\mathcal{H}$ to be \emph{mixed} when the coordinates of its centre are of opposite parities, and let also $c$ be either one of the two mixed cells of $\mathcal{H}$ with a corner at $p$.

Observe that $p$ is the midpoint of an edge of $\mathcal{U}_\text{Odd}$ in $U$ if and only if the common sides of $c$ with $c'$ and $c''$ are in the same state relative to $E$; i.e., either both belonging to $E$ or both not belonging to $E$. This is equivalent to $c'$ and $c''$ being in the same state relative to $[E]$; i.e., either both contained within $[E]$ or both outside of $[E]$.

The same reasoning applies to all pairs of one odd and one even cell of $\mathcal{H}$ which share a corner at the centre of a cell of $P$. Recall that, by virtue of $P = 2S$, we can view $P$ as the union of several $2 \times 2$ blocks corresponding to the cells of $S$. For each such block, one odd cell of $\mathcal{H}$ is concentric with it and four even cells of $\mathcal{H}$ are centred at its corners; so the odd and even cells of $\mathcal{H}$ with centres in $P$ (either on the boundary or in the interior of $P$) are all connected by their corners. Hence, $U$ is free of edges of $\mathcal{U}_\text{Odd}$ if and only if all odd cells of $\mathcal{H}$ contained within $P$ are contained within $[E]$, too; and all even cells of $\mathcal{H}$ lie outside of $[E]$.

By the analysis in Section \ref{prelim}, if $E$ is regular then these conditions are indeed satisfied. Suppose now that the conditions are satisfied; so $[E]$ is the union of all odd cells of $\mathcal{H}$ contained within $P$ together with some mixed cells of $\mathcal{H}$ contained within $P$. This is enough to guarantee that $[E]$ is of the form $\llbr 2T \rrbr$ for some spanning subgraph $T$ of the adjacency graph of $S$. Since $[E]$ is connected, so must be $T$ as well; and also $T$ must be acyclic because $[E]$ is free of holes. We conclude that $T$ is in fact a spanning tree in the adjacency graph of $S$, as desired. \end{myproof}

We call a pseudotour of $P$ \emph{turn-even} when its turn graph satisfies the conditions of Lemma~\ref{even}; i.e., when it is a subgraph of $\mathcal{U}_\text{Even}$.

%% file: mt-paper-04-down.tex
\section{The Downward Phase of the Proof} \label{down}

By Lemma \ref{even}, for our main result it suffices to show that every minimum-turn tour of $P$ is turn-even. The analogue of this claim in the pseudotour setting is quite easy to prove, and we do so below. (Notice that a ``minimum-turn tour'' is one which attains the smallest number of turns out of all tours; whereas a ``minimum-turn pseudotour'' is one which attains the smallest number of turns out of all pseudotours. The two minima are taken over distinct ranges, and will in general differ.)

\begin{proposition*} The minimum-turn pseudotours of an even polyomino are all turn-even. \end{proposition*}

\begin{myproof} Let $E$ be a pseudotour of $P$ and let $U$ be its associated turn graph. Suppose first that $U$ contains a cycle $C$. When we delete $C$ from $U$, all three conditions (i)--(iii) of Lemma \ref{tgp} are preserved. So the deletion yields a new turn graph $U'$ which corresponds to a new pseudotour $E'$ of $P$. But the number of turns in $E'$ is smaller; specifically, it has decreased by the length of~$C$. We conclude that all minimum-turn pseudotours of $P$ must be associated with acyclic turn graphs.

Suppose, then, that the turn graph $U$ of $E$ is indeed acyclic. This means that each one of its connected components is a tree. By Lemma \ref{deg}, each leaf of such a tree must be a corner of~$P$. However, since $P$ is an even polyomino, all of its corners are vertices of $\mathcal{U}_\text{Even}$. Hence, each connected component of $U$ is a subgraph of $\mathcal{U}_\text{Even}$; and so is $U$ as well. \end{myproof}

In the proof of the Proposition, we constructed an improved pseudotour directly. For the tour analogue of this result, our approach will be more roundabout. Then, as outlined in the introduction, we are going to obtain an improved tour by working our way through a series of intermediate pseudotours.

The same basic idea of deleting cycles from the turn graph in order to reduce the number of turns continues to be helpful in the tour setting, too. However, we must be a lot more careful with the deletions so as to be able to keep track of one more key parameter -- the number of cycles in the current pseudotour. We proceed now to describe one ``refined'' variant of the cycle deletion procedure which takes these considerations into account.

\begin{lemma} \label{del} Let $E$ be a pseudotour of $P$ with $s$ cycles and $t$ turns. Suppose that the turn graph of $E$ is not acyclic. Then there exists a pseudotour $E'$ of $P$ with $s'$ cycles and $t'$ turns such that $s' \le s + \ell - 1$ and $t' = t - 2\ell$ for some positive integer $\ell \ge 2$. \end{lemma}

\begin{myproof} Let $C$ be any cycle in the turn graph $U$ of $E$ such that the interior of $[C]$ does not contain edges of $U$. The existence of such a cycle is guaranteed by virtue of $P$ being free of holes.

Delete $C$ from $U$. As in the proof of the Proposition, this yields a new turn graph $U'$ associated with a new pseudotour $E'$ of $P$. Since $\mathcal{U}$ is bipartite, the length of $C$ is even -- say, $2\ell$. Clearly $t' = t - 2\ell$ because the number of turns in a pseudotour equals the number of edges in its turn graph. What remains is to estimate $s'$.

Since there are no edges of $U$ in the interior of $[C]$, all edges of $E$ contained within $[C]$ are of the same orientation -- either horizontal or vertical. Suppose, for concreteness, the former. Then the effect of the deletion of $C$ from $U$ is that these horizontal edges are all removed from $E$ and replaced with the vertical edges of $\mathcal{H}$ contained within $[C]$.

Consider any newly created cycle $D$ of $E'$, meaning any cycle which belongs to $E'$ but not to $E$. Then $D$ cannot be contained within $[C]$ because all edges of $E'$ contained within $[C]$ are vertical. But $D$ cannot lie entirely outside of $[C]$, either, because then its edges would have been unaffected by the transition from $E$ to $E'$, and it would have already been present in $E$. Hence, $D$ must include both some edges contained within $[C]$ and some edges which lie outside of $[C]$.

It follows that every newly created cycle $D$ of $E'$ must cross over $C$ at least twice. Since each edge of $C$ offers exactly one opportunity for such a crossing, we conclude that there are at most $\ell$ newly created cycles in $E'$. Furthermore, at least one cycle of $E$ must have been destroyed during the transition from $E$ to $E'$. Taken together, these observations imply that $s' \le s + \ell - 1$, as desired. \end{myproof}

By iterating the procedure of Lemma \ref{del}, one by one we can get rid of all cycles in the turn graph of $E$. Suppose that this requires $k$ iterations in total, with $\ell = \ell_1$, $\ell_2$, $\ldots$, $\ell_k$. Let also $d = \ell_1 + \ell_2 + \cdots + \ell_k$. Then, overall, $s$ will increase by at most $d - k$; whereas $t$ will decrease by exactly $2d$. Or, to summarise, we have learned how to make the turn graph of our pseudotour acyclic while exercising some degree of control over the number of turns and the number of cycles in the pseudotour itself.

%% file: mt-paper-05-up.tex
\section{The Upward Phase of the Proof} \label{up}

Our next order of business will be to devise a procedure for stitching together the cycles of pseudotours. This will be the ``upward'' counterpart of the ``downward'' procedure in Lemma~\ref{del}. Once again, we must be careful to ensure that both the number of turns and the number of cycles change in a tightly controlled manner.

\begin{lemma} \label{glue} Let $E$ be a turn-even pseudotour of $P$ with $s \ge 2$ cycles and $t$ turns. Then there exists a turn-even pseudotour $E'$ of $P$ with $s - 1$ cycles and at most $t + 2$ turns. \end{lemma}

\begin{myproof} Since the adjacency graph of $P$ is connected, it must contain two neighbouring vertices which belong to two distinct cycles of $E$. Suppose, up to rotation, that $(x, y)$ and $(x + 1, y)$ are in the cycles $C$ and $D$ of $E$, respectively.

We can assume without loss of generality that $y$ is as small as possible, over all such pairs. It follows that at least one of $(x, y - 1)$---$(x, y)$ and $(x + 1, y - 1)$---$(x + 1, y)$ is not an edge of $E$, as otherwise $(x, y - 1)$ and $(x + 1, y - 1)$ would form a different suitable pair with a smaller ordinate. Suppose, for concreteness, that $(x, y - 1)$---$(x, y)$ does not belong to $E$. Then $(x, y)$ must be joined to $(x - 1, y)$ and $(x, y + 1)$ in $E$.

Suppose, for the sake of contradiction, that $(x + 1, y)$---$(x + 1, y + 1)$ is not an edge of $E$. Then, similarly, $(x + 1, y)$ must be joined to $(x + 1, y - 1)$ and $(x + 2, y)$ in $E$. (Figure \ref{false}.) This means that both of $(x - 1/2, y + 1/2)$---$(x + 1/2, y - 1/2)$ and $(x + 1/2, y + 1/2)$---$(x + 3/2, y - 1/2)$ are edges in the turn graph $U$ of $E$. However, they belong to distinct connected components of $\mathcal{U}$ -- in contradiction with $E$ being turn-even and $U$ being a subgraph of $\mathcal{U}_\text{Even}$.

\begin{figure}[ht] \null \hfill \begin{subfigure}[c]{55pt} \centering \includegraphics{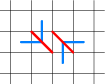} \caption{} \label{false} \end{subfigure} \hfill \begin{subfigure}[c]{105pt} \centering \includegraphics{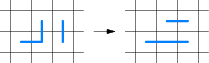} \caption{} \label{flip} \end{subfigure} \hfill \null \caption{} \label{ff} \end{figure}

We conclude that $(x + 1, y)$---$(x + 1, y + 1)$ must be an edge of $E$. Consider now the four vertices $(x, y)$, $(x, y + 1)$, $(x + 1, y)$, $(x + 1, y + 1)$ in the adjacency graph of $P$. We know that the two vertical edges $e_1$ and $e_2$ of $\mathcal{H}$ between them belong one each to the two distinct cycles $C$ and $D$ of $E$. Hence, the two horizontal edges $e_3$ and $e_4$ of $\mathcal{H}$ between them do not belong to $E$.

Delete $e_1$ and $e_2$ from $E$, and replace them with $e_3$ and $e_4$. (Figure \ref{flip}.) This \emph{flip} transforms $E$ into a new pseudotour $E'$ of $P$ where $C$ and $D$ have been stitched together into a single cycle, while the other cycles of $E$ have remained unaffected. Furthermore, at least one turn of $E$ has been substituted with a straight in $E'$ -- specifically, the one at $(x, y)$. Since our flip involves just four cells of $P$ altogether, it follows that during the transition from $E$ to $E'$ the number of turns has increased by at most $3 - 1 = 2$.

We are only left to verify that $E'$ is still turn-even. There are three potential spots where a new turn could appear during the transition from $E$ to $E'$, namely $(x, y)$---$(x + 1, y)$---$(x + 1, y - 1)$, $(x, y + 2)$---$(x, y + 1)$---$(x + 1, y + 1)$, and $(x, y + 1)$---$(x + 1, y + 1)$---$(x + 1, y + 2)$. The edges of $\mathcal{U}$ associated with them form a $4$-cycle together with the edge of $\mathcal{U}$ associated with the turn $(x - 1, y)$---$(x, y)$---$(x, y + 1)$ of $E$ which we just straightened out. Hence, all four of these edges are in the same connected component of~$\mathcal{U}$. Since $E$ is turn-even by assumption, this connected component must be $\mathcal{U}_\text{Even}$. \end{myproof}

By iterating the procedure of Lemma \ref{glue}, one by one we can stitch all cycles of $E$ together into a tour. Clearly, this is going to require $s - 1$ iterations in total. Furthermore, over the course of the entire process, $t$ will increase by at most $2s - 2$. Or, to summarise, our toolbox has come to include also a method for the transformation of turn-even pseudotours into tours; we pay for this by means of an increase in the number of turns, but -- as will become evident shortly -- the price is just right for our purposes.

Both of the ``downward'' and ``upward'' procedures that we alluded to in the introduction are now in place, and we are properly equipped to establish our main result:

\begin{theorem*} Let $P$ be an even polyomino. Then all minimum-turn tours of $P$ are regular. \end{theorem*}

\begin{myproof} Consider any irregular tour $E$ of $P$ with $t$ turns, and let its turn graph be $U$. By Lemma~\ref{even}, we get that $U$ contains an edge of $\mathcal{U}_\text{Odd}$. However, as in the proof of the Proposition, all acyclic connected components of $U$ are subgraphs of $\mathcal{U}_\text{Even}$; and so $U$ contains a cycle.

This means that we can apply the procedure of Lemma \ref{del}. We keep applying it, iteratively, until it cannot be applied any longer; i.e., until we obtain a pseudotour $E'$ of $P$ whose turn graph $U'$ is acyclic. Let $s'$ and $t'$ be the number of cycles and the number of turns in $E'$, respectively. Since the number of cycles has increased by $s' - 1$ in the transition from $E$ to $E'$, it follows that the number of turns has decreased by at least $2s'$. Hence, $t' \le t - 2s'$.

Since $U'$ is acyclic, as in the proof of the Proposition we get that it is a subgraph of $\mathcal{U}_\text{Even}$; i.e., $E'$ is turn-even. This means that we can apply the procedure of Lemma \ref{glue}. We keep applying it, iteratively, until it cannot be applied any longer; or, in other words, until we obtain a tour $E''$ of~$P$. Let $t''$ be the number of turns in $E''$. Since the number of cycles has decreased by $s' - 1$ in the transition from $E'$ to $E''$, it follows that the number of turns has increased by at most $2s' - 2$. Thus $t'' \le t' + 2s' - 2 \le t - 2$.

We conclude that the number of turns in $E''$ is strictly smaller than the number of turns in~$E$. Therefore, an irregular $E$ can never be a minimum-turn tour of $P$. \end{myproof}

It is far from obvious whether our argument can be adapted to the setting in which regions with holes are allowed. Peng and Rascoussier's conjecture remains open in full generality.